\documentclass[11pt]{article}
\newcounter{numcitacao}

  \newcommand{\bdm}{\begin{displaymath}}
 
 \newcommand{\edm}{\end{displaymath}}
  \newcommand{\be}{\begin{equation}}
  \newcommand{\la}{\langle}
  \newcommand{\ra}{\rangle}
  \newcommand{\ee}{\end{equation}}
  \newcommand{\bea}{\begin{eqnarray}}
  \newcommand{\eea}{\end{eqnarray}}
  \newcommand{\beann}{\begin{eqnarray*}}
  \newcommand{\eeann}{\end{eqnarray*}}

  \newcommand{\alg}{{\cal A}}

  \newcommand{\cp}{\begin{proof} }
  \newcommand{\tp}{\end{proof} }

\newtheorem{lemma}{Lemma}
\newtheorem{obs}[lemma]{Observation}

\newtheorem{theorem}[lemma]{Theorem}

  \newenvironment{proof}{\noindent {\bf Proof:}\ }{~\rule{1.5ex}{1.5ex}\\}

  \newcounter{numpasso}

\title{Idempotent Generated Endomorphisms of an Independence Algebra}
\author{Jo\~ao Ara\'ujo}
\date{\today}

\usepackage{amsfonts}
\usepackage{amsmath}
\usepackage{amssymb}
\usepackage[ansinew]{inputenc}

\begin{document}
\def\diag{\mathop{\rm diag}}
\def\rank{\mathop{\rm rank}}
\def\End{\mathop{\rm End}}
\def\Aut{\mathop{\rm Aut}}
\def\GL{\mathop{\rm GL}}
\def\S{\mathop{\cal S}}
\def\halmos{\hfill\rule{6pt}{6pt}}
\def\R{\mathord{\sf{l\hspace{-0.1em}R}}}
\def\C{\mathord{\hspace{0.15em}\sf{l\hspace{-0.4em}C}}}
\def\Q{\mathord{\hspace{0.15em}\sf{l\hspace{-0.4em}Q}}}
\def\Z{\mathord{\sf{Z\hspace{-0.4em}Z}}}
\def\N{\mathord{\sf{l\hspace{-0.1em}N}}}
\def\D{\mathord{\sf{l\hspace{-0.1em}D}}}
\def\P{\mathord{\sf{l\hspace{-0.1em}P}}}
\def\F{\mathord{\sf{l\hspace{-0.1em}F}}}

\pagestyle{plain}

\def\diag{\mathop{\rm diag}}
\def\rank{\mathop{\rm rank}}
\def\Ker{\mathop{\rm Ker}}
\def\End{\mathop{\rm End}}
\def\Aut{\mathop{\rm Aut}}
\def\GL{\mathop{\rm GL}}
\def\S{\mathop{\cal S}}
\def\halmos{\hfill\rule{6pt}{6pt}}
\def\R{\mathord{\sf{l\hspace{-0.1em}R}}}
\def\C{\mathord{\hspace{0.15em}\sf{l\hspace{-0.4em}C}}}
\def\Q{\mathord{\hspace{0.15em}\sf{l\hspace{-0.4em}Q}}}
\def\Z{\mathord{\sf{Z\hspace{-0.4em}Z}}}
\def\N{\mathord{\sf{l\hspace{-0.1em}N}}}
\def\D{\mathord{\sf{l\hspace{-0.1em}D}}}
\def\P{\mathord{\sf{l\hspace{-0.1em}P}}}
\def\F{\mathord{\sf{l\hspace{-0.1em}F}}}

\maketitle

\begin{abstract}

The aim of this note is to give a direct proof for the following
result proved by Fountain and Lewin: {\em Let $\alg$ be an
independence algebra of finite rank and let $a$ be a singular
endomorphism of $\alg $. Then $a=e_1\ldots e_n$ where $e^2_i=e_i$
and $rank(a)=rank(e_i)$.}

 \

{\em $2000$ Mathematics Subject Classification\/}: 20M10.
\end{abstract}

$\ $

We assume the reader to have a basic knowledge of independence
algebras. As reference we suggest \cite{gould}.  Throughout this
note $\alg$ will be an independence algebra of finite rank
($rank(\alg )\geq 2$) and with universe $A$. By $End(\alg) $ and
$PEnd(\alg )$ we will denote, respectively, the monoid of
endomorphisms and of partial endomorphisms of $\alg$, and
$Aut(\alg)$ will be the automorphism group of $\alg$. Moreover,
$End(\alg )\setminus Aut(\alg )$ is the monoid of singular
endomorphisms and will be denoted by $Sing(\alg )$. For $a\in
PEnd(\alg )$ and $X\subseteq A$, denote by $\Delta (a), \nabla
(a)$, respectively, the domain and  the image of $a$, and denote
by $a|X$ the restriction of $a$ to $X$. The following  proof is
inspired by \cite{araujosf} and goes as
 follows:
\begin{description}
  \item[(1)] for every $a\in Sing(\alg )$ there exists an
  idempotent $e\in End(\alg )$ and a partial monomorphism $a'$
  such that $a=ea'$ (with $\nabla (e)=\la E\ra=\Delta a'$ and $\nabla (a') =\la
  Ea\ra$, for some independent set $E$);
  \item[(2)] there exist idempotents $e_1, \ldots , e_n \in
  PEnd(\alg )$ such that $rank(e_i)=rank(a)$ and $(E)e_1\ldots e_n=Ea$;
  \item[(3)] for every $f\in Sym(Ea)$ (the symmetric group on $Ea$),
  there exist idempotents $f_1, \ldots , f_m \in PEnd(\alg )$ such that $rank(f_i)=rank(a)$ and $f=(f_1 \ldots
  f_m)|Ea$;
  \item[(4)] finally, using (2), $a'|E: E  \rightarrow  Ea$ can be
  factorized as $a'|E=(e_1 \ldots  e_nf)|E $, with $f\in Sym(Ea)$,
  and hence,
  by (3), we have $a'|E=(e_1 \ldots  e_nf_1 \ldots
  f_m)|E$, so that $a=ea=ea'=ee_1 \ldots  e_nf_1 \ldots
  f_m$.
\end{description}

We start by proving (1). Let $a\in Sing(\alg )$. Since $End(\alg
)$ is regular there exists an idempotent $e\in End(\alg )$ such
that $ea=a$ (and $rank(e)=rank(a)$). Moreover, if $E$ is a basis
for $\nabla (e)$, we have
\[|Ea|\leq |E|
=rank(e)=rank(a)=rank(ea)=rank(\la E \ra a)=rank(\la E a\ra )\leq
|Ea|.\]
 Therefore $rank (\la E a\ra )= |Ea|$ and hence  $Ea$ is
independent. Since $|E|=|Ea|$ and $Ea$ is independent, then  $a|E$
is $1-1$ and  $a|\la E \ra$ is a monomorphism. It is proved that
$a=e(a|\la E \ra)$ and  (1) of the scheme above follows.

For the remains of this note $a$, $e$, $E$ and $Ea$ are the
objects introduced above and are fixed. We now prove (2).   Let
$l=rank(a)$ and let $K_l =\{B\subseteq A \mid |B|=l \mbox{ and $B$
is independent}\}$. Consider in $K_l$ the following relation:
$(B_1,B_2)\in \rho$ if and only if $|B_1\setminus B_2|\leq 1$.

\begin{lemma}
For some natural number $n\geq 0$ we have $(E,Ea)\in \rho^n$.
\end{lemma}
\cp Let $\uparrow\!\! E=\{C\in K_l \mid (E,C)\in \rho^m, \mbox{
for some natural } m\}$. Now let $C\in \uparrow\!\! E$ such that
for all $D\in \uparrow\!\! E$ we have $|C\cap Ea|\geq |D\cap Ea|$.
We claim that $C=Ea$. In fact, if by contradiction $C\neq Ea$,
then there exists $d\in C\setminus Ea$ (as $|C|=|Ea|$). Therefore,
for some $c\in Ea$, we have $c\not\in \la C\setminus \{d\}\ra$
(since we cannot have $Ea\subseteq \la C\setminus \{d\} \ra$ as
$rank(\la C\setminus \{d\}\ra)<rank(\la Ea\ra )$). Thus
$C_0=(C\setminus \{d\})\cup \{c\} \in K_l$ and $(C,C_0)\in \rho$.
Since  $(E,C)\in \rho^k$ and $(C,C_0)\in \rho$ it follows that
$(E,C_0)\in \rho^{k+1}$ so that $C_0\in \uparrow\!\! E$ and
$|C_0\cap Ea|>|C\cap Ea|$, a contradiction. It is proved that
$C=Ea$ and hence $Ea\in \uparrow\!\! E$.
 \tp

We observe that the key fact used in the proof of this lemma is a
well known property of matroids (see \cite{oxley}, Exercise 1.,
p.15).

 Now let $E= E_1
\rho E_2 \rho \ldots \rho E_{n-1}\rho E_n=Ea$. We claim that, for
every $i=1,\ldots , n-1$, there exist two idempotents $e_{i,1},
e_{i,2}\in PEnd(\alg )$ such that $(E_i) e_{i,1}e_{i,2}=E_{i+1}$.
In fact let $D=E_i\cap E_{i+1}$ so that $E_i=D\cup \{x\}$ and
$E_{i+1}=D\cup \{y\}$. Then either $y\not\in \la D\cup \{x\} \ra$
or $y\in \la D\cup \{x\} \ra$.

In the first case, $D'=D\cup \{x,y\}$ is independent. Thus we can
define a mapping $f:D'\rightarrow D\cup \{y\}$ by $xf=yf=y$ and
$df=d$, for all $d\in D$. This mapping $f$ induces a partial
idempotent endomorphism $e_{i,1}:\la D'\ra \rightarrow \la D\cup
\{y\}\ra$ such that $(E_i)e_{i,1}=(D\cup \{x\})e_{i,1}=D\cup
\{y\}=E_{i+1}$. (Let $e_{i,2}=e_{i,1}$).

In the second case, $rank (\la D\cup \{x,y\} \ra )=rank(\la D\cup
\{x\} \ra )< rank (\alg )$ and hence there exists $z\not\in \la
D\cup \{x,y\} \ra $. Therefore  $D\cup \{x,z\}$ and $ D\cup
\{y,z\}$ are independent sets and hence there exist two
idempotents $e_{i,1},e_{i,2}\in PEnd(\alg )$ such that
$e_{i,1}|D=id_D=e_{i,2}|D$, $xe_{i,1}=ze_{i,1}=z$ and
$ze_{i,2}=ye_{i,2}=y$. It is obvious that $(D\cup
\{x\})e_{i,1}e_{i,2}=D\cup \{y\}$ and hence
$(E_i)e_{i,1}e_{i,2}=E_{i+1}$. We have constructed idempotents
 $e_{1,1},e_{1,2},\ldots , e_{n,1},e_{n,2}\in PEnd(\alg )$ such that
$(E)e_{1,1}e_{1,2}\ldots e_{n-1,1}e_{n-1,2}=Ea$ and hence $(2)$ is
proved.
\begin{obs}\label{obs1}
Observe that in both cases considered above, and for all $i=1,
\ldots , n-1$, we have $\nabla (e_{i,1})\leq \Delta (e_{i,2})$ and
$\nabla (e_{i,1}e_{i,2})\leq \Delta (e_{i+1,1})$.
\end{obs}
\begin{lemma}
Let $Sym(Ea)$ be the symmetric group on $Ea$ and let $(xy)\in
Sym(Ea)$ be a transposition. Then there exist idempotents $f_1,
f_2, f_3 \in PEnd(\alg )$ such that  $f=(f_1 f_2
  f_3)|Ea$ and $rank(f_i)=rank(a)$.
\end{lemma}
\cp Since $|Ea|<rank (\alg )$ there exists $z\not\in \la Ea \ra$
and hence $Ea \cup \{z\}$ is independent. Let $f_1, f_2, f_3$ be
idempotents of domain $\la Ea \cup \{z\}\ra $ defined as follows:
$xf_1=z=zf_1$ and $uf_1=u$, for the remaining elements of $Ea$;
$xf_2=x=yf_2$ and $uf_2=u$, for the remaining elements of $Ea$;
$zf_3=y=yf_3$ and $uf_1=u$, for the remaining elements of $Ea$.
Hence $(x)f_1f_2f_3=(z)f_2f_3=(z)f_3=y$ and
$(y)f_1f_2f_3=(y)f_2f_3=(x)f_3=x$. Moreover, for all $b\in
Ea\setminus \{x,y\}$, $(b) f_1f_2f_3=b$. It is proved that $(f_1
f_2 f_3)|Ea=(xy)$ and it is clear that $rank(f_i)=rank(a)$. \tp
\begin{lemma} Let
$f\in Sym(Ea)$. Then there exist idempotents $f_1, \ldots , f_m
\in PEnd(\alg )$ such that  $f=(f_1 \ldots
  f_m)|Ea$ and $rank(f_i)=rank(a)$.
\end{lemma}
\cp  Clearly  $ (Ea)f=(Ea)(x_1y_1)\ldots (x_my_m)$, since every
permutation of a finite set can be decomposed in to a product of
transpositions. Now the result follows by repeated application of
the previous lemma.
 \tp
\begin{obs}\label{obs2}
Observe that, for all the partial idempotents $f_i$ considered in
the proof of the two previous lemmas, we have $\Delta (f_i)= \la
Ea \cup \{z\}\ra$ and $\nabla (f_i)= \la Ea \ra$. Thus $\nabla
(f_i)< \nabla (f_{i+1})$.
\end{obs}
\begin{lemma}
Let $i_0,i_1 \in PEnd(\alg )$ and suppose that $\nabla (i_0) \leq
\Delta (i_1)$. Then there exist $\beta_0 , \beta_1 \in End(\alg )$
such that $(\beta_0 \beta_1 )|\Delta (i_0)=i_0i_1$, $rank(\beta_j
)=rank(e_j)$ and if $i_j $ is idempotent, then $\beta_j$ is
idempotent ($j\in \{0,1\}$).
\end{lemma}
\cp Let $i_0:\la B_0\ra \rightarrow \la C_0 \ra$ be a partial
endomorphism ($B_0, C_0$ are independent sets). Then we can extend
$B_0, C_0$, respectively, to $B,C$, bases of $\alg$ and define
$\beta_0:\alg \rightarrow \alg$ such that $\beta_0|B_0=i_0$ and
$(B\setminus B_0)\beta_0=\{c\}\subseteq C_0$. In the same way we
extend $i_1$ to $\beta_1\in End(\alg )$.  Therefore
$\beta_0|\Delta (i_0)=i_0$ and $\nabla (\beta_0)=\nabla (i_0)$.
Thus $(\beta_0i_1)|\Delta(i_0)=i_0i_1$ and
  $(\beta_0 \beta_1 )|\Delta (i_0)=i_0i_1$.  \tp

\begin{theorem} (Fountain and Lewin \cite{fou2})
Every $a\in Sing(\alg )$ is the product of idempotents $e_1,\ldots
,e_n \in End(\alg )$ such that $rank(a)=rank(e_i)$.
\end{theorem}
\cp Let $a\in Sing (\alg )$. Then there exists an idempotent $e$
(with $\nabla (e)=\la E \ra$) and a partial monomorphism $a'=a|\la
E\ra$ such that $a=ea'$. Moreover $a'$ maps the basis $E$ of
$\nabla (e)$ into a basis $Ea$ of $\nabla (a)$. We proved that for
some idempotents $e_1, \ldots , e_k\in PEnd(\alg )$ we have
$(E)e_1 \ldots e_k=Ea$. Let $h=e_1 \ldots e_k$ and consider $f\in
Sym(Ea)$ defined by $(xh)f=xa$, for all $x\in E$. It is obvious
that $(e_1 \ldots e_kf)|E=a|E$ and hence $(e_1 \ldots e_k\phi)|\la
E\ra =a|\la E\ra =a'$, where $\phi \in Aut(\la Ea \ra )$ is (the
automorphism) induced by $f$. Therefore $ee_1 \ldots
e_k\phi=ea'=a$. It is also proved that for some idempotents $f_1,
\ldots , f_n \in PEnd(\alg )$ we have $f=(f_1 \ldots f_n)|Ea$.
Thus $\phi=(f_1 \ldots f_n)|\la Ea\ra$ and hence
\[a=ea'=ee_1
\ldots e_k\phi=ee_1 \ldots e_k(f_1 \ldots f_n)|\la Ea\ra=ee_1
\ldots e_kf_1 \ldots f_n
\] (since $\nabla (ee_1 \ldots e_k)=\la Ea\ra$).
 To finish the proof we only need to show
that the idempotent partial endomorphisms used above can be
replaced by idempotent total endomorphisms.  This follows from
Observations \ref{obs1} and \ref{obs2} and the previous lemma.\tp

 {\em Acknowledgments: I wish to express my deepest gratitude to
 the referee for suggesting many improvements to this note and also to
 Professor Gracinda Gomes.}

\
\begin{center}
\begin{tabular}{ll}
   Universidade Aberta &\hspace{4cm}  Centro de Álgebra \\
  R. Escola Politécnica, 147 &\hspace{4cm} Universidade de Lisboa\\
  1269-001 Lisboa Portugal &\hspace{4cm} 1649-003 Lisboa Portugal\\
    & \hspace{4cm} mjoao@lmc.fc.ul.pt
\end{tabular}
\end{center}


\begin{thebibliography}{99}
\bibitem{araujosf} J.\ Araújo, {On Idempotent Generated Semigroups}, Semigroup Forum, {\bf 65} (2002), 138-140.



\bibitem{fou2}
J.\ Fountain, A.\ Lewin, {\em Products of idempotent endomorphisms
of an independence algebra of finite rank},  Proc.\ Edinburgh\
Math.\ Soc.\ {\bf 35} (1992), 493--500.


\bibitem{gould}
V.\ A.\ R.\ Gould, {\em Independence algebras}, Algebra Universalis {\bf 33} (1995), 327--329.

\bibitem{oxley}
J.\ G.\ Oxley, {\em Matroid Theory}, Oxford University Press,
(1992).
\end{thebibliography}
\end{document}